\documentclass[12pt]{article} 
\usepackage{amssymb}

\usepackage{latexsym}
\usepackage{epsfig}

\newif{\ifcomentarios}
\comentariosfalse

\renewcommand{\Bbb}{\mathbb}

\begin{document}
\title{Rakib--Sivashinsky and Michelson--Sivashinsky Equations for
Upward Propagating Flames: A Comparison Analysis } 
\author{ \textbf{Leonardo F. Guidi}\thanks{ 
Supported by FAPESP under grant 98/10745-1. E-mail: {\it
guidi@if.usp.br}} 
\ and \ \textbf{Domingos H. U. Marchetti}\thanks{
Partialy supported by CNPq and FAPESP. E-mail:{\it \
marchett@if.usp.br}} \\ 
Instituto de F\'{\i}sica\\
Universidade de S\~{a}o Paulo\\
Caixa Postal 66 318\\
05315 S\~{a}o Paulo, SP, Brazil}
\date{\today }
\maketitle

\begin{abstract}
We establish a comparison between Rakib--Sivashinsky and
Michelson--Sivashinsky quasilinear parabolic differential equations governing
the weak thermal limit of upward flame front propagating in a channel. For
the former equation, we give a complete description of all steady solutions
and present their local and global stability analysis. For the latter,
multi--coalescent unstable steady solutions are introduced and shown to be
exponentially more numerous than the previous known coalescente solutions.
This fact is argued to be responsible for the disagreement of the observed
dynamics in numerical experiments with the exact (linear) stability analysis
and also gives the ingredients to describe the quasi--stable behavior of
parabolic steadily propagating flame with centered tip.
\end{abstract}

\newpage

Rakib--Sivashinsky quasilinear parabolic differential equation governing the
weak thermal limit of an upward flame interface propagating in a channel is
given by 
\begin{equation}
\phi _{t}=\varepsilon \phi _{xx}-\frac{1}{2}\phi _{x}^{2}+\phi -\stackrel{\_
}{\phi }\,,  \label{RSeq}
\end{equation}
where $y=\phi (t,x)$, $0<x<\pi $, defines an instantaneous flame profile in
dimensionless variables, $\stackrel{\_}{\phi }=
{\displaystyle{1 \over \pi }} \displaystyle\int 
_{0}^{\pi }\phi (t,x)\,dx$ denotes the space average, $\varepsilon >0$ is a
physical parameter (Markstein length) and Neumann (adiabatic) boundary
condition at channel walls is imposed: $\phi _{x}(t,0)=\phi _{x}(t,\pi )=0$.
According to Darrieus--Landau's hydrodynamic flame theory\cite{RS} there are
two competing sources of instabilities in this model given by the last two
terms of (\ref{RSeq}). A plane flame front separating the cold and hot gas
is subject to Rayleigh--Taylor instability, due to the thermal expansion,
and to the buoyancy effect caused by external acceleration. These two
ingredients lead the front to become convex toward the cold gas. As in \cite
{RS}, the acceleration vector points toward positive $y$--values and the
upward propagating direction is negative. For the purpose of comparison such
convention will be maintained.

Paraboloid profiles with the tip located around the center of the channel
are observed, both experimentally and in computational simulations, to
remain ``stable'' for long time (see e.g. \cite{MiS} and references
therein). There 
are, however, experiments whose paraboloid tip eventually slides to the
channel wall. A mathematical description of the former configuration as a
``quasi--equilibrium transient state'' was provided by Berestycki, Kamin and
Sivashinsky \cite{BKS} in which work stationary solutions of (\ref{RSeq}),
their respective stability properties and the nonlinear dynamics were
described for sufficiently small $\varepsilon $. Metastable dynamics has
been studied in details in \cite{SW}. The proofs of the other theorems
stated in \cite{BKS} remain, to the best of our knowledge, unpublished.

Sivashinsky'{}s previous equation of weak thermal expansion\cite{S,MS} 
\begin{equation}
\phi _{t}=\varepsilon \phi _{xx}-\frac{1}{2}\phi _{x}^{2}+I(\phi )\,,
\label{MSeq}
\end{equation}
hasn'{}t considered buoyancy effect and Darrieus--Landau instability
has been taken 
into account replacing $\phi -\stackrel{\_}{\phi }$ in (\ref{RSeq}) by a
linear singular integral operator\cite{nota0} given by a multiplication by $
\left| k\right| $ in the Fourier representation 
\[
I\left( \phi \right) (t,x)=\frac{1}{\pi }\sum_{k=1}^{\infty }k\,\int_{-\pi
}^{\pi }\cos \left[ k\left( x-y\right) \right] \,\phi (t,y)\,dy\,. 
\]
In contradistinction, the only trace of $I(\phi )$ term in equation (\ref
{RSeq}) comes from the removal of $k=0$ mode 
\[
\phi -\stackrel{\_}{\phi }=\frac{1}{\pi }\sum_{k=1}^{\infty }\int_{-\pi
}^{\pi }\cos \left[ k\left( x-y\right) \right] \,\phi (t,y)\,dy\,\,. 
\]
Michelson--Sivashinsky equation (\ref{MSeq}) has been studied by many
authors after Thual--Frish--H\'{e}non's application\cite{TFH} of pole
decomposition. For a more recent survey of this method, see Vaymblat and
Matalon\cite{VM1,VM2} and references therein.

Despite of the fact that equations (\ref{RSeq}) and (\ref{MSeq}) differ in
many respects, \ it is our purpose to expose the similarities and
distinctions of their solutions. We shall see, by the introduction of so
called bi--coalescent steady states, that equation (\ref{MSeq}) may also
admit quasi--equilibrium transient ``parabolic'' profile with centered tip.
Our comparison relies on two recent analysis. Firstly, Vaymblat and Matalon 
\cite{VM1} have determined all steady solutions whose poles coalesce into a
line parallel to the imaginary axis and solved explicitly the eigenvalue
problem of (\ref{MSeq}) linearized about those states. They conclude, in a
second paper\cite{VM2}, the existence of only one linearly stable
steady coalescent pole solution.

For equation (\ref{RSeq}) with a different parametrization, we have given
\cite{GM} a complete description of all equilibrium solutions and
provided their local and global stability in an appropriated Sobolev space.
Since our results are major extensions of those stated in \cite{BKS} and the
mathematical presentation in \cite{GM} may cause certain difficulties in
translating to the present application, we shall here restate them with a
brief explanations of their proofs. A detailed presentation including the
metastable states analysis will appear elsewhere \cite{MG}.

We shall first restrict ourselves to the Rakib--Sivashinsky equation. The
details presented after each statement (detached by brackets) are
essentials to establish the subsequent comparison but may be skipped
in a first reading. Multi--coalescent steady states of
Michelson--Sivashinsky equation will be considered next. A conclusion
will be presented at the end. 

To discuss our results let us consider the equation 
\begin{equation}
u_{t}=\varepsilon u_{xx}-u\,u_{x}+u\,\,,  \label{u}
\end{equation}
for the derivative $u=\phi _{x}$, with Dirichlet boundary conditions $
u(t,0)=u(t,\pi )=0$. There is a one--to--one correspondence between 
\begin{equation}
\theta (t,x):=\int_{0}^{x}u(t,y)\,dy\,,  \label{phiu}
\end{equation}
with $u$ a solution of (\ref{u}), and\ a solution $\phi $ of (\ref{RSeq})
given as follows. Notice that $\theta $ satisfy (\ref{RSeq}) with $\stackrel{
\_}{\phi }$ replaced by $\varepsilon \phi _{xx}(t,0)$. By definition, $
\theta (t,0)=0$ for all $t\geq 0$ and their steady solutions do not
propagate (they are equilibrium). We have $\theta (t,x)=\phi (t,x)-\phi
(t,0)\,$, $\theta _{x}=\phi _{x}$ and $\theta _{xx}=\phi _{xx}$. For the
opposite relation, we have 
\[
\phi +\int_{0}^{t}\stackrel{\_}{\phi }(s)\,e^{t-s}\,ds=\theta +\varepsilon
\int_{0}^{t}\theta _{xx}(s,0)\,e^{t-s}\,ds\,, 
\]
because both sides satisfy equation (\ref{RSeq}) without $\stackrel{\_}{\phi 
}$ term, which yields 
\begin{equation}
\phi (t,x)=\theta (t,x)-\int_{0}^{t}\left( \stackrel{\_}{\theta }(\tau
)\,-\varepsilon \,\theta _{xx}(\tau ,0)\right) \,\,ds\,.
\label{phitheta} 
\end{equation}
This and (\ref{phiu}) will be used to discuss the steadily propagating
solutions of (\ref{RSeq}).

If $A$ denotes the operator given by r.h.s. of (\ref{u}) linearized about
the trivial solution $u_{0}=0$, let ${\cal B}^{1/2}$ be the Banach space
equipped with the graph norm $\left\| f\right\| _{1/2}=
\displaystyle\int 
_{0}^{\pi }\left| A^{1/2}f(x)\right| ^{2}\,dx$ (equivalent to the Sobolev
space $H_{0}^{1}(0,\pi)$). In \cite{GM}, Theorems 3.2 and 5.10, we have proven:

\noindent 
{\it The initial value problem (\ref{u}) with }$u(0,\cdot )=u_{0}\in
{\cal B}^{1/2}$
{\it \ has a unique solution for all }$t>0${\it \ and the trajectories }$
\left\{ u(t,\cdot )\right\} _{t\geq 0}${\it \ lie in a compact set in }$
{\cal B}^{1/2}${\it .}

By (\ref{phiu}) and (\ref{phitheta}), this provides the existence and
uniqueness of the initial value problem (\ref{RSeq}) for all $t>0$.
Compactness property will be useful for the global stability analysis.

The steadily propagating solutions of (\ref{RSeq}) can be obtained from the
equilibrium solution of (\ref{u}). The quantity inside the parenthesis
in (\ref{phitheta}) remains positive under the dynamics of (\ref{u})
and the 
front propagates upward (toward negative $y$--values). As a consequence, the
flame front propagates steadily with velocity 
\begin{equation}
V=-\,\stackrel{\_}{\vartheta }\,+\varepsilon \,\vartheta _{xx}(0)=-\, 
\frac{1}{2} \; \overline{\vartheta ^{2}},  \label{vv} 
\end{equation}
provided $\vartheta (x)=\int_{0}^{x}v(y)\,dy$, where $v$ is a non trivial
solution of 
\begin{equation}
\varepsilon v_{xx}-v\,v_{x}+v=0\,,  \label{v}
\end{equation}
with $v(0)=v(\pi )=0$. The next result is Theorem 4.1 of \cite{GM}:

\noindent 
{\it For }$\varepsilon >1${\it , }$v_{0}=0${\it \ is the unique
solution of (\ref{v}). For }$\varepsilon <1${\it \ such that }$ 
{\displaystyle{1 \over k+1}}
\leq \varepsilon <
{\displaystyle{1 \over k}}
${\it \ holds for some }$k\in N${\it , there exist }$2k${\it \ nontrivial
solutions }$u_{j}^{\pm }${\it , }$j=1,\ldots ,k${\it , with }$j-1${\it \
zeros in }$\left( 0,\pi \right) ${\it \ and such that }$%
u_{j}^{-}(x)=-u_{j}^{+}(\pi -x)${\it \ for }$j${\it \ odd and }$%
u_{j}^{-}(x)=u_{j}^{+}(\pi /j+x)${\it \ mod }$\pi ${\it \ for }$j${\it \
even. Each pair }$u_{j}^{\pm }${\it \ bifurcate from the trivial solution }$%
v_{0}=0${\it \ at }$\varepsilon _{j}=1/j^{2}${\it \ with }$%
\lim\limits_{\varepsilon \uparrow \varepsilon _{j}}u_{j}^{\pm }=0${\it . In
the phase space }${\Bbb R}^{2}${\it , }$\left( u_{j}^{\pm \prime
},u_{j}^{\pm }\right) ${\it \ describes closed orbits around }$\left(
0,0\right) ${\it \ whose distance from the origin increases monotonically as 
}$\varepsilon ${\it \ decreases.}

\begin{figure}[!ht] 
\begin{center}
\epsfig{file=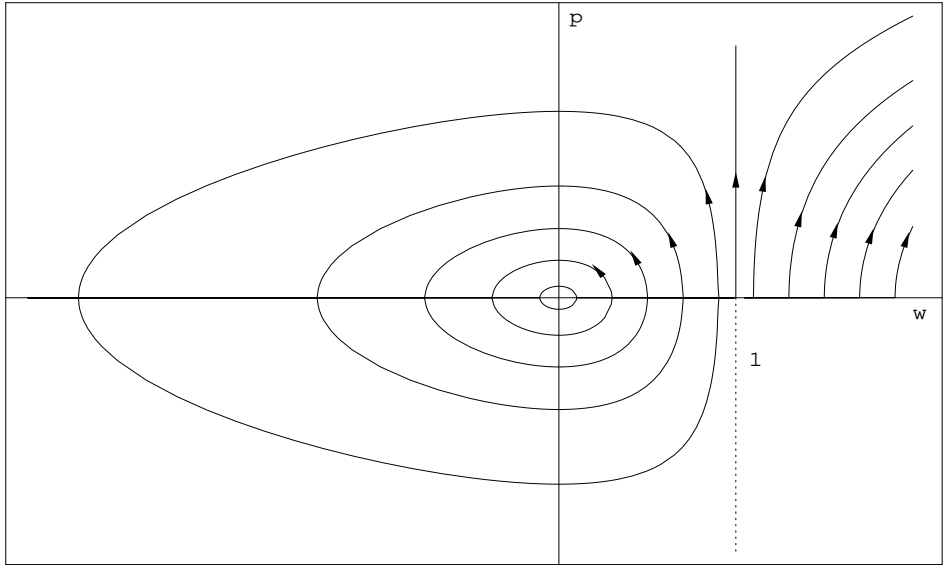, height=7.5cm, width=12.5cm}
\end{center}
\caption{Trajectories of the dynamical system (\ref{ds}).} 
\label{phsp} 
\end{figure}

\noindent [Equation (\ref{v}) can be written as a dynamical
system  
\begin{equation}
\left\{ 
\begin{array}{lll}
w \acute{}
& = & \varepsilon ^{-1}p\left( w-1\right)  \\ 
p \acute{}  = & w\,,
\end{array}
\right.   \label{ds}
\end{equation}
where $p=v$ and $w=v{\acute{}}$. Since (\ref{ds}) remains unaltered by
changing $p\rightarrow -p$ and $ 
x\rightarrow -x$, the orbits are symmetric with respect to the $w$--axis and
Dirichlet boundary conditions $p(0)=p(\pi )=0$ hold for any periodic orbit
with period $2\pi $ such that $p=0$ at $x=0$. In \cite{GM} we have shown
that {\bf (i)} ${\Bbb R}^{2}$ is foliated by non--overlapping orbits 
\[
\gamma _{w_{0}}=\left\{ (w(x),p(x)):x\in {\Bbb R}\;{\rm and}
\;(w(0),p(0))=\left( w_{0},0\right) \right\} 
\]
labeled by the coordinate $w_{0}$ of the positive $w$--axis; {\bf (ii)} as $
\varepsilon $ varies, the orbit varies continuously from one to another; 
{\bf (iii) }the orbits are closed if $w_{0}<1$ and open if $w_{0}\geq 1$; 
{\bf (iv)} the origin $\left( 0,0\right) $ is enclosed by any closed orbits.
The trajectories $\gamma _{w_{0}}$, $w_{0}\geq 0$, are portrayed in
Figure \ref{phsp}.

The period $T=T(\varepsilon ,w_{0})$ of a closed orbit varies continuously
in $\varepsilon $ and $w_{0}$ and is given by 
\[
T=\int_{\gamma _{w_{0}}}dx=2\int \frac{dp}{w}\,, 
\]
where the second integral is calculated in the semi--orbit. Hence, $T=2\pi $
defines implicitly a continuous function $w_{0}=w_{0}(\varepsilon )$, whose
corresponding semi--orbit $\gamma _{w_{0}}=\left\{ \left( w(x),p(x)\right)
:0\leq x\leq \pi \right\} $ determines the solution $v(x)=p(x)$ of
(\ref{v}). The question to be addressed now is how many $2\pi
$--periodic orbits 
there are for each $\varepsilon $. If $w_{0}$ is small, $\gamma _{w_{0}}$
can be approximated by an ellipsoid: $\varepsilon ^{-1}p^{2}+w^{2}=w_{0}^{2}$
; 
\begin{equation}
T_{{\rm ellipse}}=4\int_{0}^{\sqrt{\varepsilon }}\frac{dp}{\sqrt{
1-\varepsilon ^{-1}p^{2}}}=2\pi \sqrt{\varepsilon }  \label{Te}
\end{equation}
uniformly in $w_{0}$ and we have 
\begin{equation}
\lim\limits_{w_{0}\rightarrow 0}T(\varepsilon ,w_{0})=2\pi \sqrt{\varepsilon 
}\,.  \label{Tw0}
\end{equation}
In \cite{GM}, we have shown 
\begin{equation}
\frac{\partial T}{\partial w_{0}}>0\,  \label{delT}
\end{equation}
holds form any $\varepsilon >0$. Consequently, $T=T(\varepsilon ,w_{0})$ is
a (strictly) monotonic increasing function of $w_{0}$ and $\varepsilon $,
since $T(\varepsilon ,w_{0})=\sqrt{\varepsilon }T(1,w_{0})$. As $
T$ increases when $\varepsilon $ and $w_{0}$ increase, equation (\ref{Tw0})
implies that there is no nontrivial solutions ($v\neq 0$) for $\varepsilon
\geq 1$, since the period of each nontrivial orbit exceeds $2\pi $. Now, let 
\[
G_{1}(\varepsilon ,w_{0})=T(\varepsilon ,w_{0})-2\pi \,, 
\]
be defined for $0<\varepsilon \leq 1$ and $0\leq w_{0}<1$. Under the
condition (\ref{delT}), there exist a {\bf unique} function $g_{1}:\left[ 0,1
\right] \longrightarrow {\Bbb R}_{+}$ with $g_{1}(1)=0$ such that $%
G_{1}(\varepsilon ,g(\varepsilon ))=0$ for all $0<\varepsilon \leq 1$. In
view of (\ref{delT}), $g(\varepsilon )$ is a monotonically (strictly)
decreasing function and $\lim_{\varepsilon \rightarrow 0}g(\varepsilon )=1$.

Now, $v\equiv 0$ is the unique equilibrium solution of (\ref{u}) for any $
\varepsilon \geq 1$. For $\varepsilon <1$ such that $1/(k+1)^{2}\leq
\varepsilon <1/k^{2}$ holds for some integer $k\geq 1$, one can apply the
implicit function theorem to equation 
\[
G_{j}(\varepsilon ,w_{0}):=T(\varepsilon ,w_{0})-\frac{2\pi }{j}=0\,,
\]
with $j=1,\ldots ,k$ and conclude, exactly as in the case for $G_{1}$, the
existence of a unique monotone decreasing function $g_{j}:\left[ 0,1/j^{2}
\right] \longrightarrow {\Bbb R}_{+}$ with $g_{j}(1/j^{2})=0$ and $
\lim\limits_{\varepsilon \rightarrow 0}g_{j}(\varepsilon )=1$, such that $
G_{j}(\varepsilon ,g_{j}(\varepsilon ))=0$. The branch $g_{j}(\varepsilon )$
defines a nontrivial equilibrium solution $v_{j}$ of (\ref{u}) given by the $
p$--component of the orbit $\gamma _{g_{j}(\varepsilon )}$ with $\left(
g_{j}(\varepsilon ),0\right) $ at $x=0$ and $p(x)=0$ at each semi--period $
x=m\pi /j$, $m=1,\ldots ,j$. Hence, if $\varepsilon $ is such that $
1/(k+1)^{2}\leq \varepsilon <1/k^{2}$, there exist $k$ equilibrium (not
identically $0$) solutions of (\ref{u}), $v_{1}^{+},\,\ldots ,$ $v_{k}^{+}$,
with $v_{j}^{+}$ having $j-1$ zeros in the interval $(0,\pi )$. Notice that $
g_{j}(\varepsilon )=v_{j}^{\prime }(0)>0$, $j=1,\ldots ,k$, and $v_{1}(x)>0$
for all $0<x<\pi $. These will be useful for the stability analysis. There
is an additional set of $k$ equilibrium solution $v_{1}^{-},\,\ldots
,v_{k}^{-}$, with $v_{j}^{-}(x)=p(\pi /j+x)$.]

When $\varepsilon $ is small enough the orbit spends most of time at right
semi--plane, $p(x)\simeq x$ mod $2\pi /j$ and $\vartheta _{j}^{\pm
}(x):=\int_{0}^{x}v_{j}^{\pm }(y)\,dy\,$\ has\ wrinkles of parabolic shape
separated by ``cusps''. From (\ref{vv}) the steadily propagating solutions
of (\ref{RSeq}) has the profile of $\vartheta _{j}$ propagating with the
velocity $V$: 
\begin{equation}
\varphi _{j}(x)=\vartheta
_{j}(x)-\,
\frac{t}{2} \; \overline{\vartheta 
_{j}^{2}}\,.  \label{phij} 
\end{equation}

The plane flame front solution $\varphi _{0}={\rm Const}$ is an
asymptotically stable solution if $\varepsilon \geq 1$ for both
equations (\ref{RSeq}) and (\ref{MSeq}). Its local stability can be
read directly from 
the spectrum $\sigma (A)$ of the operator $A$, given by (\ref{RSeq})
and (\ref{MSeq}) linearized about $\varphi _{0}$, which consist of
simple eigenvalues 
\begin{equation}
\lambda _{n}=1-\varepsilon n^{2}\qquad {\rm and}\qquad \eta
_{n}=n-\varepsilon n^{2},\qquad n=1,2\ldots \,,\,  \label{eigenv}
\end{equation}
respectively, with corresponding eigenfunction $\sqrt{1/\pi }\cos nx$. So, $
\lambda _{n}<0$ and $\eta _{n}<0$ if $\varepsilon >1$. At $\varepsilon =1$,
the trivial solution $\varphi _{0}$ bifurcates \cite{nota} into a steadily
propagating front in which, when extended periodically to the real line as
an even function, there is one ``parabolic'' tip centered at one wall and a
single ``cusp'' at another. The next result of \cite{GM}, Theorems $5.1$ and 
$5.14$, states that such configuration is the only one globally
asymptotically stable solution of (\ref{RSeq}) for all $\varepsilon <1$.
Parabolic front with centered tip or wrinkled flame fronts may be described
as quasi--equilibrium states discussed in refs. \cite{MS,BKS,SW}. By (\ref
{phitheta}), it is enough to examine the stability of the nontrivial
equilibrium solutions of (\ref{u}).

\noindent 
{\it If }$u_{0}\in {\cal B}^{1/2}${\it \ and }$\varepsilon >1${\it , then the
initial value problem (\ref{u}) with }$u(0,\cdot )=u_{0}${\it \ satisfies }$
\lim\limits_{t\rightarrow \infty }\left\| u(t,\cdot )\right\| _{1/2}=0$ and
the trivial solution is globally stable{\it . For }$\varepsilon <1${\it \ so
that }$1/(k+1)^{2}\leq \varepsilon <1/k^{2}${\it \ holds, there exist }$\rho
>0${\it \ such that, if }$\left\| u_{0}-v\right\| _{1/2}\leq \rho ${\it ,
then }$\lim\limits_{t\rightarrow \infty }\left\| u(t,\cdot )-v\right\|
_{1/2}=0${\it \ for }$v=v_{1}^{\pm }${\it \ and, for any sequence }$\left\{
u_{n}\right\} _{n\geq 1}${\it \ with }$\lim\limits_{n\rightarrow \infty
}\left\| u_{n}-v\right\| _{1/2}=0${\it , we have }$\sup\limits_{t>0}\left\|
u(t,\cdot )-v\right\| _{1/2}\geq \delta >0${\it \ for all }$n${\it \ and }$v$
{\it \ equal to }$0${\it \ or any equilibrium }$v_{j}^{\pm }${\it , }$
j=2,\ldots ,k${\it . Moreover, there is an open dense set }$U\subset
{\cal B}^{1/2}$ 
{\it \ containing the origin and }$\left\{ v_{j}^{\pm }\right\} _{j=1}^{k}$
{\it such that, if }$u_{0}\in U${\it \ and }$\varepsilon <1${\it , then }$
\lim\limits_{t\rightarrow \infty }\left\| u(t,\cdot )-v_{1}^{\pm }\right\|
_{1/2}=0${\it .}

\noindent [Let us first examine the linear stability. 
If $u(t,x)=v(x)+\zeta (t,x)$ then equation (\ref{u}) can be
written as 
\[
\zeta _{t}=L\zeta +\zeta \zeta _{x}\,,
\]
where\cite{nota3} 
\[
L\xi =L\left[ v\right] \xi =\varepsilon \xi ^{\prime \prime }-v\xi ^{\prime
}+(1-v^{\prime })\xi 
\]
is the linearization of the right hand side of (\ref{u}) about $v$. Acting
on the space of functions $\xi $ satisfying $\xi (0)=\xi (\pi )=0$, $L$ is
symmetric, $\left( \eta ,L\xi \right) _{\rho }=\left( L\eta ,\xi \right)
_{\rho }$, with respect to the inner product 
\[
\left( \eta ,\xi \right) _{p}=
\displaystyle\int 
_{0}^{\pi }\eta (x)\,\xi (x)\,\rho (x)\,dx
\]
with weight $\rho (x)=\exp \left\{ -\varepsilon
^{-1}\int_{0}^{x}v(y)\,dy\right\} $. \ As a consequence, we can apply the
comparison theorem to establish the following criterium (see \cite{GM} for
details).

\noindent Stability Criterium. {\it If} $\varphi $ {\it is the solution of} 
\[
L\left[ v\right] \varphi =0\,, 
\]
{\it on} $0<x<\pi $, {\it satisfying} $\varphi (0)=0$ {\it and} $\varphi
^{\prime }(0)=1$, {\it the largest eigenvalue} $\lambda $ {\it of }$L\left[ v
\right] $ {\it is negative if} $\varphi >0$ {\it on} $\left( 0,\pi \right) $ 
{\it and positive if there exist} $x^{\prime }$ {\it such that }$0<x^{\prime
}<\pi $ {\it and} $\varphi (x^{\prime })=0$.

For an equilibrium solution $v$ of (\ref{u}), let 
\[
\chi =c\left( 2v-v^{\prime \prime }\right) \,,  
\]
where $c>0$ is chosen so that $\chi ^{\prime }(0)=1$. It follows from
equilibrium equation $\varepsilon v^{\prime \prime }=-\left( 1-v^{\prime
}\right) v$ that $\chi (0)=0$ and $\chi >0$ whenever $v>0$ (recall $v(0)=0$
and $\left( 1-v^{\prime }\right) >0$ for all closed orbits). Moreover, if $
v>0$, an explicit calculation gives (see \cite{GM}) 
\[
L\left[ v\right] \chi =-2\left( v^{\prime }\right) ^{2}v<0\,\,,  
\]
and this implies, by applying the comparison theorem once more, 
\[
\varphi (x)\geq \chi (x) 
\]
for all $0\leq x<\pi $. It thus follows from the stability criterium that $
v_{1}$ is a stable equilibrium solution for all $\varepsilon <1$ (recall $
v_{1}(x)>0$ for all $0<x<\pi $, implying $\chi (x)>0$ in the same domain).
To show $v_{2},\ldots ,v_{k}$ are unstable, we observe
\begin{equation}\label{Lv}
L\left[ v_{j}\right] v_{j}^{\prime } =\left( \varepsilon v_{j}^{\prime
\prime }+(1-v_{j}^{\prime 
})v_{j}\right) ^{\prime }=0\,, 
\end{equation}
by the equilibrium equation. It thus follows that $v_{j}$ with $j\geq 2$ is
unstable in view of the equilibrium criterium and the fact that $v_{j}$ has
at least one zero in $\left( 0,\pi \right) $ and the same holds for $\varphi 
$ by the Wronskian positivity $W(\varphi ,v_{j}^{\prime };x)>0$. It is
important to note that $v_{j}^{\prime }$ fails to be an eigenfunction of $L
\left[ v_{j}\right] $ with $0$ eigenvalue because it does not satisfies the
boundary conditions required. Approximate eigenfunction with exponentially
small (in $\varepsilon $) eigeinvalue can, however, be constructed using
boundary layer techniques (see e. g. \cite{SW}). Equation (\ref{Lv}) holds
also for any stationary solution $v$ of Michelson--Sivashinsky equation.

We turn to the global stability. A Liapunov function for equation (\ref{u})
can be constructed via the generalized Euler--Lagrange method due to
Zelenyak, Lavrentiev and Vishnevskii \cite{ZLV}. Let 
\[
U(u)=\int_{0}^{\pi }\Phi (u,u_{x})\,dx 
\]
where $\Phi (p,w)=-p^{2}/(2\varepsilon )+\left( 1-w\right) \ln \left(
1-w\right) +w$ is an appropriated ``Lagrangian'' . The total derivative of $
U $ with respect to $t$ is obtained by the calculus of variation 
\begin{eqnarray*}
\stackrel{\cdot }{U}\left( u\right) &=&-\int_{0}^{\pi }\left(  
\frac{\partial}{\partial x} \frac{\partial \Phi }{\partial
u_{x}}-\frac{\partial \Phi v
}{\partial u} \right) u_{t}\,dx \\   
&=&-\int_{0}^{\pi }\rho (u_{x})\left( \varepsilon  
u_{xx}-u\,u_{x}+u\right) u_{t}\,dx  
\end{eqnarray*}
where $\rho (w)=\varepsilon ^{-1}/(1-w)$ is a positive weight. Note $
\stackrel{\cdot }{U}$ is negative in view of (\ref{u}) and $U$ is a Liapunov
functional. Due to the fact that the trajectories $\left\{ u(t,\cdot
)\right\} _{t\geq 0}$ lie in a compact set, LaSalle's invariance principle
can be applied to show that all solutions of (\ref{u}) in ${\cal B}^{1/2}$
converge to an equilibrium solution as $t\rightarrow \infty $. From
equations (\ref{phitheta}), (\ref{phiu}) and the two results stated before
one concludes that $\varphi _{0}=0$ is globally asymptotically stable
solution of (\ref{RSeq}) for $\varepsilon \geq 1$. In addition, if $
\varepsilon <1$, there is an open dense set ${\cal U}\subset {\cal B}^{1/2}$
of initial condition such that all solutions $\phi (t,x)$, with initial
condition on ${\cal U}$, is asymptotically of the form $\vartheta _{1}+Vt$
where $\vartheta _{1}=\int_{0}^{x}v_{1}(y)\,dy$, the velocity $V=\stackrel{\_\_
}{-\vartheta _{1}^{2}}/2$ and $v_{1}$ is the stable equilibrium
solution of (\ref{u}).]

\medskip

\begin{center}
\noindent MULTI--COALESCENT\ POLE\ SOLUTIONS
\end{center}

\medskip

The same scenario seems to hold for equation (\ref{MSeq}), at least if one
restricts to the space of coalescent pole solutions. According to \cite
{TFH,VM1}, if $\varepsilon $ is such that $%
{\displaystyle{1 \over 2n+1}}
\leq \varepsilon <
{\displaystyle{1 \over 2n-1}}
$, $n=1,2,\ldots $, there exist$\ n$ steadily propagating fronts, $\varphi
_{j}(t,x)=\vartheta _{j}(x)+Vt$, with $\vartheta _{j}(x)=
\displaystyle\int 
_{0}^{x}v_{j}(y)\,dy$ and 
\[
v_{j}(x)=\varepsilon \sum_{k=1}^{2j}\cot \left( \frac{x-z_{k}}{2}\right) \,
\]
for $j=1,\ldots ,n$, each corresponding to the coalescent $j$--pole solution
of (\ref{MSeq}), but only the one with largest number of poles, $\varphi _{n}
$, is asymptotically stable. When $\varepsilon $ crosses $1/(2n+1)$ from
above, $\varphi _{n}$ becomes unstable and a new solution $\varphi _{n+1}$
bifurcates from the former solution.

$\varphi _{j}$ was called coalescent $j$--pole solution because of the
dynamic of poles $\left\{ z_{k}(t)\right\} _{k=1}^{2j}$, induced by (\ref
{MSeq}), tends to align them parallel to the imaginary axis. In \cite
{TFH,VM1}, (\ref{MSeq}) was considered with periodic boundary conditions.
With Newmann (adiabatic) boundary conditions there are two sets of steadily
propagating fronts $\left\{ \varphi _{j}^{\pm }\right\} _{j=1}^{n}$
distinguished by the location of their poles: $\varphi _{j}^{+}$ and $
\varphi _{j}^{-}$ has poles aligned at $\Re e\,z=0$ and $\Re e\,z=\pi $,
respectively.

Note that there is a relation between the number of steadily propagating
fronts and of positive eigenvalues in the spectrum $\sigma (A)=\left\{ \eta
_{j}\right\} _{j\geq 1}$ of the linear operator $Av=\varepsilon v^{\prime
\prime }+I(v)$, since the latter increases by one unit at $\varepsilon
_{n}=1/n$ and the former at $\widetilde{\varepsilon }_{n}=1/\left(
2n-1\right) $. For equation (\ref{RSeq}), the number of steady solutions and
the dimension of the unstable manifold ${\cal M}_{n}$ of the trivial
solution $\varphi _{0}={\rm const}$, are equal with the instability value
now located at $1/n^{2}$. There, in contradistinction, $\varphi _{1}^{\pm }$
are the only asymptotically stable solution for all $\varepsilon <1$.

Although (\ref{phij}) holds for $\varphi
_{j}^{\pm }$ with $v_{j}^{\pm }$ a coalescent pole solution of 
\begin{equation}
\varepsilon v_{xx}-v\,v_{x}+I(v)=0\,,  \label{vI}
\end{equation}
$v_{j}^{\pm }$ does not vanish in $\left( 0,\pi \right) $. One may define a
coalescent $j$--pole solution $v_{j}^{(k)\pm }$ with $k-1$ zeros in $\left(
0,\pi \right) $, $k=1,2,\ldots $, by setting $v_{j}^{(k)\pm }(x)=
{\displaystyle{1 \over k}}
v_{j}^{\pm }(kx)$ mod $\pi $. Note $v_{j}^{(k)\pm }$ solves (\ref{vI})
provided $v_{j}^{\pm }$ solves the same equation with $\varepsilon $
replaced by $k\varepsilon $. As a consequence, a sets of steadily
propagating fronts $\left\{ \varphi _{j}^{(k)\pm }\right\} _{j=1}^{n}$
exists if $\varepsilon $ is such that $
{\displaystyle{1 \over 2n+1}}
\leq k\varepsilon <
{\displaystyle{1 \over 2n-1}}
$ holds for some $n=1,2,\ldots $ (see Figure \ref{bif}). In total,
there are $2 \mathop{\displaystyle\sum} _{m=1}^{n}\left[ {2n / (2m-1)}
\right] = o(n^{2})$ coalescent steady solutions for $\varepsilon \gtrsim 
{\displaystyle{1 \over 2n+1}}
$, where $\left[ z\right] $ means the integer part of real number $z$. The
solutions $\varphi _{j}^{(k)\pm }$ with $k\geq 2$, are not stable and may
represent the cellular profile observed experimentally provided an associate
quasi--equilibrium solution described in \cite{MS,BKS,SW} can be defined.

\begin{figure}[!ht] 
\begin{center}
\epsfig{file=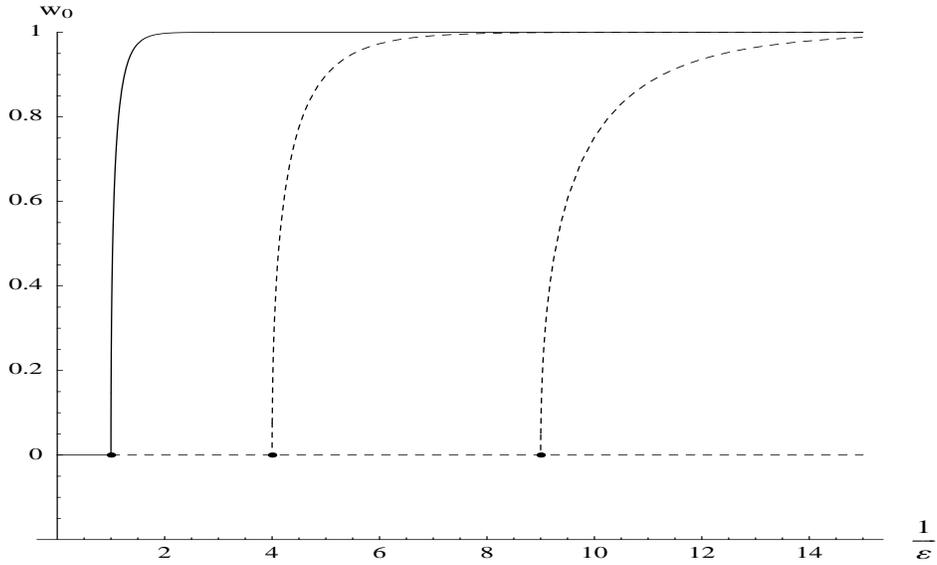, height=7.5cm, width=12.5cm}
\end{center}
\caption{Bifurcation diagram for the distance $\Delta \varphi $
from cusp to tip of steadily state solutions of (\ref{RSeq}). Solid
and dashed lines refer, respectively, to stable 
and unstable solutions.}
\label{bifur} 
\end{figure}

\begin{figure}[!ht] 
\begin{center}
\epsfig{file=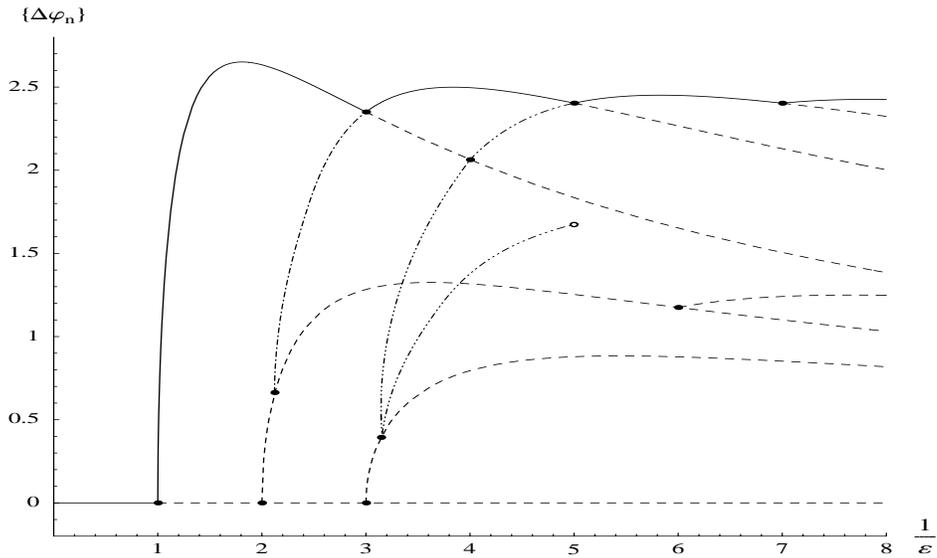, height=7.5cm, width=12.5cm}
\end{center}
\caption{Bifurcation diagram for the distance $\Delta \varphi $
from cusp to tip of steadily state solutions of (\ref{MSeq}). Solid
and dashed lines refer, respectively, to stable 
and unstable coalescente solutions. Doted-dashed lines refer to
unstable multi--coalescente solutions.} 
\label{bif} 
\end{figure}

In the following, for each $\varepsilon $ such that $
{\displaystyle{1 \over n+1}}
\leq \varepsilon <
{\displaystyle{1 \over n}}
$ holds, a new family ${\cal F}_{n}$ of steadily propagating flame front
solutions, denominated multi--coalescent $n$--pole solutions, will be
introduced. Our preliminary investigation indicates that there are at
least an exponential number $c^{n-1}$ of solutions in ${\cal F}_{n}$ and all,
but $2$ of them, seems to be 
unstable by numerical computation. As a consequence, the invariant set $
{\cal K}_{n}^{MS}=\bigcup_{\varphi }{\cal M}_{n}(\varphi )$, defined as the
union of the unstable manifold of all equilibrium solutions, for
equation (\ref{MSeq}) differs enormously from the invariant set ${\cal
K}_{n}^{RS}=
{\cal M}_{n}(\varphi _{0})$ for equation (\ref{RSeq}). Here, $n$ indicates
the number of bifurcations with respect to the trivial solution $\varphi _{0}
$. In particular, ${\cal K}_{n}^{MS}$ may have dimension exponentially more
numerous than the dimension of ${\cal K}_{n}^{RS}$ (for comparison,
see Figures \ref{bifur} and \ref{bif}). We believe that this crucial
distinction is responsible for the disagreement between the numerical
study by Gutman--Sivashinky \cite{GS} and the exact calculation by
Vaymblat--Matalon \cite{VM1}. 

The bi--coalescent $n$--pole solutions $\varphi _{n_{0},n_{\pi }}$ are
indexed by $(n_{0},n_{\pi })$ with $n_{0}+n_{\pi }=n$ indicating the number
of pairs of complex conjugate poles\cite{nota4} in each line $\Re e\,z=0$
and $\Re e\,z=\pi $. Note that the system of equations governing the
dynamics of the poles $z_{j}=x_{j}+iy_{j}$, $j=1,\ldots ,2n$, in the complex
plane, given by 
\[
\stackrel{\cdot }{z_{j}}=-\varepsilon \sum_{l\neq j}\cot \left( \frac{
z_{j}-z_{l}}{2}\right) -i\,\frac{y_{j}}{\left| y_{j}\right| }\,,
\]
preserves the location of real part $x_{j}$ since, in this case, $\Re
e\,\cot \left( 
{\displaystyle{z_{j}-z_{l} \over 2}}
\right) =0$. The poles of the bi--coalescent solution $\varphi
_{n_{0},n_{\pi }}$ thus satisfies 
\begin{equation}
\stackrel{\cdot }{x_{j}}=0\qquad {\rm and}\qquad \stackrel{\cdot }{y_{j}}
=F_{j}  \label{xy}
\end{equation}
where 
\begin{equation}
F_{j}=\varepsilon \sum_{l\neq j}\left( \coth \left( \frac{y_{j}-y_{l}}{2}
\right) \right) ^{\eta _{jl}}-\frac{y_{j}}{\left| y_{j}\right| }\,.
\label{Fj}
\end{equation}
with $\eta _{jl}:=\cos \left( x_{j}-x_{l}\right) $ taking $\pm 1$ values
according the poles $z_{j}$ and $z_{l}$ are in the same or different line.
Because of the real parte $x_{j}$ of the pole gives the ``cusp'' position of
a propagating flame profile, the bi--coalescent solution has its tip
centered somewhere in the interior of channel.

Thual, Frisch and H\'{e}non \cite{TFH} have proven that, provided $n$ is
such that $\varepsilon \left( 2n-1\right) <1$, there exist one and only one
coalescent steady solution and any solution of (\ref{xy}) with $\eta _{jl}=1$
for all $j$ and $\,l$, tends toward this steady state as $t\rightarrow
\infty $. The proof of these properties is based on the existence of a
Liapunov function with negative curvature in every direction. Here, there
exists a Liapunov function for bi--coalescent solutions 
\[
U=\varepsilon \sum_{j,l:j\neq l}\ln \left| \frac{e^{(y_{j}-y_{l})/2}+\eta
_{jl}\,e^{-(y_{j}-y_{l})/2}}{2}\right| -\sum_{j}\left| y_{j}\right| \,\,,
\]
satisfying $\stackrel{\cdot }{U}=
\mathop{\displaystyle\sum}
_{j=1}^{2n}F_{j}^{2}\geq 0$ whose Hessian matrix $
H=[U_{y_{i}y_{j}}]_{i,j=1}^{2n}$ cannot be \ proven to be negative definite
in the case of strictly bi--coalescent solution ($n_{0},\,n_{\pi }\neq 0,n/2$
) since its Ger\v{s}gorin discs may have non--vanishing intersection with
the semi--plane $\Re e\,z\geq 0$.  As $U$ may have several local maxima and
saddle points in this case (see Figure 2), any solution of equations (\ref
{xy}) tends toward to a steady bi--coalescent state as $t\rightarrow \infty $
but uniqueness cannot be guaranteed.

Let the poles $z_{1},\ldots ,z_{2n}$ of a bi--coalescente steady solution $
\varphi _{n_{0},n_{\pi }}$ be indexed as follows: $x_{j}=0$ if $j=1,\ldots
,n_{0}$, $x_{j}=\pi $ if $j=n_{0}+1,\ldots ,n$ and $y_{n+j}=-y_{j}$ for $
j=1,\ldots ,n$. The case with $n_{0}=n_{\pi }=n/2$ plays special role to
describe the stability of coalescent solutions. Note that, if $%
y_{j}=y_{j+n/2}$ for $j=1,\ldots ,n/2$, then 
\begin{eqnarray*}
v_{n/2,n/2}(x) &=&\varepsilon \sum_{j=1}^{n/2}\sum_{\eta \in \left\{
-1,1\right\} }\left\{ \cot \left( \frac{x-i\eta y_{j}}{2}\right) +\cot
\left( \frac{x-\pi -i\eta y_{j}}{2}\right) \right\}  \\
&=&2\varepsilon \sum_{j=1}^{n/2}\sum_{\eta \in \left\{ -1,1\right\} }\cot
\left( x-i\eta y_{j}\right) 
\end{eqnarray*}
corresponds to a coalescent $n/2$--pole solution $v_{n/2}^{2}$ with one zero
in $\left( 0,\pi \right) $. As $\varepsilon $ varies from $1/\left(
2n-1\right) $ to $1/\left( 2n+1\right) $, a point $\left\{ z_{j}\right\}
_{j=1}^{2n}$ satisfying $y_{j}=y_{j+n/2}$ can be shown to change from saddle
point to a global maximum of $U$, turning a local maximum somewhere in
between.

Two conclusions can immediately be extracted from these observations.
Because (\ref{xy}) tends to align the poles along either the line $\Re e\,z=0
$ or $\Re e\,z=\pi $, the coalescent $n$--pole solution $\varphi _{n}$ is
more stable than $\varphi _{n/2}^{2}$ since $\varphi _{n/2}^{2}$ is
equivalent to a bi--coalescent $n$--pole solution $\varphi _{n/2,n/2}$ which
is unstable under small perturbation which involves the real part $x_{j}$ of
the poles. Moreover, we may construct from a coalescent $n$--pole solution $
\varphi _{n}^{k}$ with $k-1$ zeros in $\left( 0,\pi \right) $ a family
of bi--coalescent solutions $\varphi _{n_{0},n_{\pi 
/k}}^{k}$ with $n_{0}$ poles aligned in $\Re e\,z=0$ and $n-n_{0}$ poles
aligned in $\Re e\,z=\pi /k$ which agree with the coalescent solution at $
\varepsilon \gtrsim 1/\left( 2kn+k\right) $ if $n_{0}=$ $n_{\pi /k}=n/2$ and 
$y_{j}=y_{j+n/2}$. Proceeding in a similar fashion, one can introduce
muti--coalescent $n$--pole solutions $\varphi _{n_{0},n_{2\pi
/k},\ldots n_{2(k-1)\pi /k} }$, with $n_{2j\pi /k} \ge 1$ poles at $\Re
e\,z=2j\pi /k$ and $
\mathop{\displaystyle\sum}
_{j=0}^{k-1}n_{2j\pi /k}=n$ and which vanishes at $k-1$ points in $\left(
0,\pi \right) $. Hence, a trajectory $\left\{ z_{j}(t)\right\} _{t\geq 0}$
of (\ref{xy}) in the phase space ${\Bbb C}^{n}$, with $\left\{
z_{j}(0)\right\} $ close to the poles of a multi--coalescent steady solution
may go along many intermediate steady states before it reaches the final
equilibrium.

\medskip 

\begin{center}
\noindent CONCLUSIONS
\end{center}

\medskip

Whether the cellular structure (flame profile with many wrinkles), observed 
numerically by Gutman--Sivashinky \cite{GS} and experimentally by \cite{G},
could be produced by the instability of the coalescent pole solutions for $
\varepsilon $ small, has been debated in the
literature (see e.g.\cite{KOP,RAS}). The work of Vaymblat--Matalon
\cite{VM1,VM2} has resolved the controversies by proving that there
always exist 
a unique (linearly) stable coalescent pole solution for $\varepsilon <1$. In 
\cite{VM1}, the discrepancy between the numerical and the exact results is
explained as an artefact of truncation and we shall not observe different
profiles if more modes were included. According to Joulin \cite{J}, once
equation (\ref{MSeq}) is incapable to describe observed wrinkled propagating
flames (as an equilibrium solution) it should be replaced by another model.

A different scenario has been presented for Rakib--Sivashinsky equation.
Numerical integrations of (\ref{RSeq}) have agreed with the analytic
prediction since the beginning, although truncated equation has been used.
Besides, parabolic profile with centered tip and cellular profiles can be
successfully explained as metastable states \cite{MiS,BKS,SW}. So, the
question to be addressed is why equation (\ref{MSeq}) is more sensitive than
(\ref{RSeq}) to be treated numerically and whether the cellular profiles can
be described as a quasi--equilibrium solution of (\ref{MSeq}).

Based in the present analysis of equation (\ref{RSeq}) and in the existence
of multi--coalescent steady states of (\ref{MSeq}), we argue in the following
that many questions remain to be investigated before Michelson--Sivashinsky
equation is abandoned.

Using the analysis in \cite{GM}, global existence in a Sobolev space $
H_{0}^{1}$ (for all times $t>0$) and uniqueness can be established for
equation (\ref{MSeq}). A basic question is: {\it Does the solution} $\phi $ 
{\it of (\ref{MSeq}) with initial condition }$\phi _{0}${\it \ in a dense
subset of} $H_{0}^{1}$ {\it converge, as} $t\rightarrow \infty ${\it \ and
for all} $\varepsilon <1$, {\it to one of coalescent steady state described
in \cite{VM1}}?

To answer this question a geometric analysis, as given for equation 
(\ref{RSeq}), must be carried out for equation (\ref{MSeq}). A
family ${\cal F}_{n}$ of steady multi--coalescent poles solutions have been
described. It may be difficult to determine whether ${\cal F}_{n}$ exhausts
all steady solutions of (\ref{MSeq}) but it is already a remarkable
difference between both equations. If all states of ${\cal F}_{n}$, except
the states $\varphi _{n}^{\pm }$ described by \cite{VM1}, are shown to be
unstable, then the unstable manifold ${\cal K}_{n}^{MS}$ containing the
origin would have dimension at least exponentially large with $n$. 

The method of \cite{ZLV} may be useful to construct a Liapunov function. To
establish the existence of a dense subspace of initial conditions from which
the solutions of (\ref{MSeq}) converge to a steady solution requires,
besides a Liapunov function, that the trajectories remain in a compact set.
For this, it is enough that the trajectories remain bounded, which follows
if the maximum principle can be shown to be adapted for equation (\ref{MSeq}).

Finally, the existence of multi--coalescent pole solutions would explain the
discrepancy between the numerical integration by Sivashinky and the linear
stability analysis in \cite{VM1}. For Rakib--Sivashinsky equation, there
exist metastable solutions whose time interval they remain ``stable'' becomes
exponentially long when $\varepsilon $ is small, creating the illusion that
they have reached the equilibrium. If metastable states can be constructed
from the bi--coalescent states is a question to be investigate. It would, in
particular, describe the quasi--stable behavior of parabolic steadily
propagating flame with centered tip. Moreover, the effect of truncation
would become more sensitive than for equation (\ref{RSeq}) in view of the
fact there is exponentially more numerous (meta)states available.


\begin{thebibliography}{BKS}
\bibitem{nota0}  $I(\phi )={\rm P.V}.
\displaystyle\int 
{\displaystyle{dy \over x-y}}
\,\phi _{x}(t,y)$ is the composition of minus the Hilbert transform and the
derivative of $\phi $ with respect to $x$.
\bibitem{nota}  In \cite{VM1}, the Darrieus--Landau instability term $I(\phi
)$ of (\ref{MSeq}) has a $1/2$ factor in front so, the bifurcation point is
at $\varepsilon =1/\gamma =1/2$. In \cite{BKS}, equation (\ref{RSeq}) is
defined for $x\in \left( 0,1\right) $ and the bifurcation occurs at $
\varepsilon =\pi ^{-2}$.
\bibitem{nota3}  Here, $L$ is defined to be minus the $L$ operator in \cite
{GM}.
\bibitem{nota4}  Because we are looking for real valued solutions of (\ref
{MSeq}), poles must occur in pairs of complex conjugate.
\end{thebibliography}
\end{document}